\documentclass[a4paper,leqno,draft]{article}

\usepackage[latin1]{inputenc}
\usepackage{latexsym}
\providecommand\ufootnote[1]{{\let\thefootnote\relax\footnote[0]{#1}}}

\newcommand{\pointir}{\discretionary{.}{}{.\kern.35em---\kern.7em}}

\newcommand{\eg}{\emph{e.g.\,}}

\newcommand{\R}{\mbox{I\hspace{-0.15em}R}}
\newcommand{\N}{\mbox{I\hspace{-0.15em}N}}
\newcommand{\Z}{\mbox{Z\hspace{-0.3em}Z}}

\newcommand{\HH}{\mbox{I\hspace{-0.15em}H}}

\newcommand{\cB}{{\mathcal B}}
\newcommand{\cC}{{\mathcal C}}
\newcommand{\cD}{{\mathcal D}}
\newcommand{\cE}{{\mathcal E}}

\newcommand{\cR}{{\mathcal R}}

\newcommand{\ip}[2]{\langle #1 , #2\rangle}
\newcommand{\aco}{\lbrace}
\newcommand{\acf}{\rbrace}

\newcommand{\pf}{\par{\noindent\textbf{Proof.~}}}

\newcommand{\qedd}{\hfill\hbox{\hskip 1pt\vrule width 4pt height 4pt
            depth 1.5pt\hskip 1pt}}

\newtheorem{theorem}{Theorem}[section]

\newtheorem{proposition}{Proposition}[section]
\newtheorem{lemma}{Lemma}[section]
\newtheorem{definition}{Definition}[section]

\newcommand{\ind}{\mbox{Ind}}

\newcommand{\Nind}{\mbox{Ind}}
\newcommand{\NDind}{\mbox{Ind}}
\newbox\qedbox
\setbox\qedbox=\hbox{$\Box$}

\title{Index growth of hypersurfaces with constant mean curvature}

\author{Pierre Bérard, Levi Lopes de Lima, Wayne Rossman}

\date{}

\begin{document}
\maketitle

\ufootnote{\hskip-0.55cm \emph{Math. Subject Classification
2000}~: 53A10, 53A35. \newline \emph{Keywords}: Constant mean
curvature, Morse index.}

\ufootnote{\hskip-0.55cm \emph{An expanded version of this paper
is available at the first and last authors' web pages.}}

\thispagestyle{empty}

\vspace{-0.5in}

\begin{abstract}

\noindent \textbf{Abstract.~} In this paper we give the precise
index growth for the embedded hypersurfaces of revolution with
constant mean curvature (cmc) $1$ in $\R^{n}$ (Delaunay
unduloids). When $n=3$, using the asymptotics result of Korevaar,
Kusner and Solomon, we derive an explicit asymptotic index growth
rate for finite topology cmc $1$ surfaces with properly embedded
ends. Similar results are obtained for hypersurfaces with cmc
bigger than $1$ in hyperbolic space.

\medskip

\noindent \textbf{Résumé.~} Dans cet article, nous estimons de
manière précise la croissance de l'indice des hypersurfaces de
révolution plongées, de courbure moyenne constante (cmc) égale à
$1$, dans $\R^{n}$ (onduloïdes de Delaunay). Quand $n=3$,
utilisant le résultat de Korevaar, Kusner et Solomon, nous en
déduisons une estimée de la croissance de l'indice des surfaces de
topologie finie, de cmc $1$ et dont les bouts sont proprement
plongés. Nous obtenons des résultats similaires pour les
hypersurfaces à cmc strictement plus grande que $1$ dans l'espace
hyperbolique.

\end{abstract}


\section{Introduction}\label{S:intro}

A complete cmc nonminimal surface without boundary in  $\R^3$ has
finite index if and only if it is compact \cite{LR}, \cite{S}. If
it is noncompact, the index is infinite, so it is natural to ask
at what rate the index grows to infinity on an exhaustion of the
surface by bounded regions. In this paper, we prove, under some
natural geometric conditions, that certain complete non-compact
cmc hypersurfaces have linear index growth. Let us give a typical
statement:

\medskip

Let $M \subset \R^3$ be a complete properly embedded
finite-topology cmc-$1$ surface. The finitely many ends $E_j, j=1
\ldots N$, of $M$ were shown by \cite{KKS} to be asymptotic to
Delaunay unduloids $\cD(\mu_j)$ with weight parameters $\mu_j
> 0$. Let $T(\mu_j)$ denote the period of the Delaunay unduloid
$\cD(\mu_j)$ and let $B(R)$ be the radius $R$ ball in $\R^3$
centered at the origin.

\begin{theorem}\label{T:thm1a} With $M$ as above, the asymptotic growth of the
index of $M$ is given by

\begin{equation}
\lim_{R\to \infty}\frac{\ind(M\cap B(R))}{R} =
\sum_{j=1}^{N}\frac{2}{T(\mu_j)}\,.
\end{equation}
\end{theorem}

\smallskip

Using the fact that $2 \le T(\mu_j) \le \pi$,
we can conclude from the preceding theorem that the
index growth provides upper and lower bounds on the number of ends of
the surface.

\smallskip

There are many known surfaces to which this theorem (or Theorem
\ref{T:thm2}) applies. Complete finite-topology cmc-$1$ surfaces
with asympotically Delaunay ends have been constructed by
N.~Kapouleas \cite{K}, R.~Mazzeo et al. \cite{MP}, and
K.~Grosse-Brauckmann et al. \cite{gks}. And there are other works
in progress for constructing such surfaces (\eg that of
J.~Dorfmeister, H.~Wu, I.~McIntosh, M.~Kilian and N.~Schmitt).
Furthermore, the structure of such surfaces is well understood
\cite{KK}, \cite{KKS}.

\medskip

The rough idea of the proof is to decompose the surface into
components, one which is a fixed compact part and the others which
are compact pieces of ends and are close to parts of Delaunay
unduloids, and then to apply Dirichlet-Neumann bracketing. We need to
show that the indexes of these end pieces
are close to the indexes of the actual Delaunay pieces, and then the heart
of the proof becomes to carefully study the indexes of the Delaunay
pieces (with both Dirichlet and Neumann boundary conditions).

\medskip

\textbf{Remark.~} Dirichlet--Neumann bracketing can be
applied to other situations. We can for example prove quadratic or cubic
index growth for certain infinite-topology cmc surfaces in
\cite{K} (Subsection 5.4).

\medskip

In Section \ref{S:frame}, we
describe the framework of the paper. In Section \ref{S:delau} we
recall the basic facts on Delaunay unduloids (in Euclidean and
hyperbolic space) and we define some special domains on
them. Section \ref{S:iecpdh} is devoted to
estimating the index of these special domains.
The main results are stated in
Subsection
5.3, and the other subsections of Section \ref{S:igr}
contain technical results needed in the proofs.

\section{Framework}\label{S:frame}

We consider hypersurfaces $M^n$ with cmc $H$ in the simply connected
$(n+1)$-dimensional space
form $\overline{M}^{\,n+1}$ with constant sectional curvature $c\in
\aco -1,0\acf$. We assume $H> |c|$. Such
hypersurfaces are critical for a variational problem whose
associated second order \emph{stability
operator} is

\begin{equation}\label{E:so1}
L := \Delta - nc - ||B||^2 \; ,
\end{equation}
where $||B||$ is the norm of the second fundamental form of the
immersion, and $\Delta$ is the (non-negative) Laplace--Beltrami
operator for the induced metric on $M$.  (When
$\overline{M} = \R^3$, we have $L = \Delta + 2K - 4H^2$,
where $K$ is the Gauss curvature.)

\smallskip

For $M$ \emph{compact}, we define $\ind(M)$ as the index (number of
negative eigenvalues) of the
quadratic form $\int_M u \, Lu \, dv_M$ on some subspace
$$
\aco u \in H^1(M) ~\big|~ u|_\Gamma = 0 \acf \; ,
$$
where $\Gamma \subset \partial M$ is a portion of the boundary of
$M$ (this means that we consider the operator $L$ with Dirichlet
boundary conditions on $\Gamma$ and with Neumann boundary conditions on
$\partial M \setminus \Gamma$). The choice of $\Gamma$ will be clear from the
context.

\smallskip

For $M$ non-compact, $\ind(M)$ is defined as the supremum of
$\ind(\Omega)$ over all relatively compact subregions $\Omega
\subset M$ (for a fixed choice of $\Gamma$).

\medskip

\textbf{Remark:~} We will not need to take into
account that there are actually two different notions of
index for cmc hypersurfaces (see \cite{BB}, \cite{BdCE},
\cite{LiRo}). Indeed, for compact subsets these indexes differ by
at most one, so their asymptotic properties are the same.

\section{Delaunay unduloids}\label{S:delau}

Here we describe Delaunay unduloids
with nonzero cmc in Euclidean and hyperbolic space.

\paragraph{3.1 Delaunay unduloids in Euclidean space, with cmc
$H > 0$.}\label{SS:dhes}
Consider a rotation hypersurface $M$ in $\R^{n+1}$ parametrized by

\begin{equation}\label{E:rs1}
\R \times S^{n-1} \ni (x,\omega) \to F(x,\omega) = \big( x, f(x)
\, \omega\big) \; .
\end{equation}

We assume $f > 0$ and $f$ is defined on $(-\infty ,
\infty )$. We choose the unit normal vector as

\begin{equation}\label{E:rsnv}
N(x,\omega) = \big( 1+ f'^2(x)\big) ^{-1/2} \, \big( f'(x), - \omega
\big) \; .
\end{equation}

\medskip

Assume that $f' \not \equiv 0$ and fix the normalized mean
curvature to be $H=1$. The profile curves of Delaunay unduloids
are given by the differential equation

\begin{equation}\label{E:delau2a}
\mu = \frac{f^{n-1}(x)}{(1 + f'^2(x))^{1/2}} - f^n(x) \,.
\end{equation}

where $\mu \in \big( 0, \frac{1}{n}(\frac{n-1}{n})^{n-1}\big)$.
The extreme values correspond to a chain of spherical beads of
radii 1 (when $\mu=0$), and to a cylinder with radius $n-1 \over
n$ (when $\mu = \frac{1}{n}(\frac{n-1}{n})^{n-1}$).

\medskip

Given $\mu \in \big( 0, \frac{1}{n}(\frac{n-1}{n})^{n-1}\big)$,
let $a_±(\mu)$ be the two positive roots of the equation $X^n -
X^{n-1} + \mu = 0$ with $a_-(\mu) \leq a_+(\mu)$.

\medskip

Let $\cD(\mu)$ be the Delaunay unduloid with cmc
$1$ and weight para\-meter $\mu \in (0,
\frac{1}{n}(\frac{n-1}{n})^{n-1}]$, whose profile curve $f$
satisfies Equation (\ref{E:delau2a}). One can show that the
function $f$ is defined over $\R$, pinched between the two
positive values $a_±(\mu)$,

\begin{equation}\label{E:delau5}
a_-(\mu) \le f(x) \le a_+(\mu),
\end{equation}

and $T(\mu)$-periodic, where $T(\mu)$ is the distance between two
consecutive values of $x$ at which $f$ achieves its least value
$a_-(\mu)$. (For a true cylinder,
$\mu=\frac{1}{n}(\frac{n-1}{n})^{n-1}$ and $f(x) = a_-(\mu) =
a_+(\mu) = \frac{n-1}{n}$ is constant. In that case,
$T(\frac{1}{n}(\frac{n-1}{n})^{n-1})$ is the limiting value of
$T(\mu)$ as $\mu$ increases up to
$\frac{1}{n}(\frac{n-1}{n})^{n-1}$.)

\medskip

The stability operator of the $n$-dimensional Euclidean Delaunay
unduloid $\cD(\mu)$ is given by

\begin{equation}\label{E:del-e5}
        L = \Delta - V \; , \;\;\;\;\;\;\;\;\;\; \mbox{where} \;\;
        V = \|B\|^2 = n \left( 1 + (n-1) \mu^2 \,f^{-2n}
        \right)\;.
\end{equation}

\begin{lemma}\label{L:new1}
  For $n\ge 2$ and for any $x\in \R$ the function $V$ in equation
(\ref{E:del-e5}) satisfies $V(x) \, f^2(x) \le n^2$.
\end{lemma}

\pf We have already seen that the weight parameter $\mu$ of
$\cD(\mu)$ satisfies $0 < \mu \le \frac{1}{n}
(\frac{n-1}{n})^{n-1}$. Consider the polynomial $P(t) = t^n -
t^{n-1} + \mu$, whose positive roots are the numbers $a_±(\mu)$.
The function $P(t)$, considered on the domain $\R_+$, achieves its
non-positive minimum $\mu - \frac{1}{n} (\frac{n-1}{n})^{n-1}$ at
$t = \frac{n-1}{n}$. Since $P(\mu^{1/(n-1)})$ and $P(1)$ are both
positive, it follows that

$$
\mu^{1/(n-1)} \le a_-(\mu) \le \frac{n-1}{n} \le a_+(\mu) \le 1\,.
$$

Consider the function $Q(t) = n t^2 (1 + (n-1) \mu^2 t^{-2n})$,
for $t>0$. When $t$ varies from $0$ to $\infty$, $Q$ decreases
from $\infty$ to its minimum $Q((n-1)^{1 \over n} \mu^{1 \over n}) \ge
0$ and then increases to $\infty$. It follows immediately that, for all
$x\in \R$,

$$
(V f^2)(x) \le \max \aco Q(a_-), Q(a_+)\acf \le \max \aco
Q(\mu^{1/(n-1)}), Q(1)\acf\,.
$$

Using the fact that $\mu \le \frac{1}{n} (\frac{n-1}{n})^{n-1}$,
it follows that $V f^2 \le n^2$ on $\R$ as claimed. \qedd


\bigskip


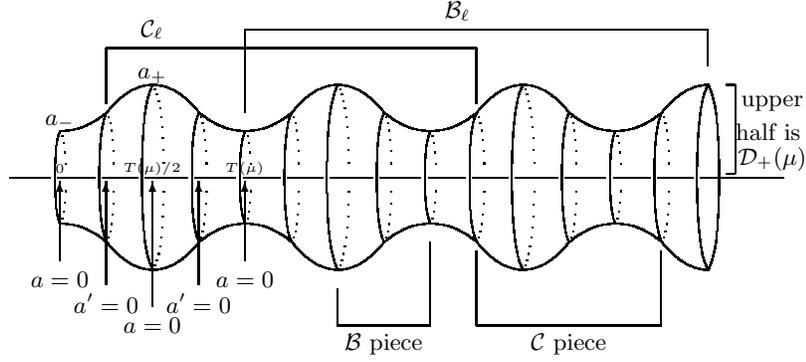
\begin{figure}[ht]
\begin{center}
\unitlength=0.7pt
\vspace{-1in}
\begin{picture}(-200.00,200.00)(0.00,300.00)
\put(-275.00,245.00){\makebox(0,0)[cc]{\small $a=0$}}
\put(-250.00,232.00){\makebox(0,0)[cc]{\small $a^\prime=0$}}
\put(-225.00,220.00){\makebox(0,0)[cc]{\small $a=0$}}
\put(-200.00,232.00){\makebox(0,0)[cc]{\small $a^\prime=0$}}
\put(-175.00,245.00){\makebox(0,0)[cc]{\small $a=0$}}
\put(-275.00,305.00){\makebox(0,0)[cc]{\tiny $0$}}
\put(-225.00,305.00){\makebox(0,0)[cc]{\tiny $T(\mu)/2$}}
\put(-175.00,305.00){\makebox(0,0)[cc]{\tiny $T(\mu)$}}
\put(-275.00,330.00){\makebox(0,0)[cc]{\small $a_-$}}
\put(-225.00,355.00){\makebox(0,0)[cc]{\small $a_+$}}
\put(-225.00,380.00){\makebox(0,0)[cc]{\small ${\cal C}_\ell$}}
\put(-60.00,390.00){\makebox(0,0)[cc]{\small ${\cal B}_\ell$}}
\put(-100.00,210.00){\makebox(0,0)[cc]{\small ${\cal B}$ piece}}
\put(0.00,210.00){\makebox(0,0)[cc]{\small ${\cal C}$ piece}}
\put(110.00,340.00){\makebox(0,0)[cc]{\small upper}}
\put(110.00,325.00){\makebox(0,0)[cc]{\small half is }}
\put(110.00,310.00){\makebox(0,0)[cc]{\small ${\cal D}_+(\mu)$}}
\put(-275,255){\vector(0,1){43}}
\put(-250,242){\vector(0,1){56}}
\put(-225,230){\vector(0,1){68}}
\put(-200,242){\vector(0,1){56}}
\put(-175,255){\vector(0,1){43}}
\put(-250,340){\line(0,1){30}}
\put(-50,340){\line(0,1){30}}
\put(-250,370){\line(1,0){200}}
\put(-175,335){\line(0,1){45}}
\put(75,360){\line(0,1){20}}
\put(-175,380){\line(1,0){250}}
\put(-125,240){\line(0,-1){20}}
\put(-75,265){\line(0,-1){45}}
\put(-125,220){\line(1,0){50}}
\put(-50,260){\line(0,-1){40}}
\put(50,260){\line(0,-1){40}}
\put(-50,220){\line(1,0){100}}
\put(85,350){\line(1,0){5}}
\put(85,302){\line(1,0){5}}
\put(90,302){\line(0,1){48}}
\bezier100(-175,325)(-162,325)(-150,335)
\bezier100(-125,350)(-138,350)(-150,335)
\bezier100(-175,275)(-162,275)(-150,265)
\bezier100(-125,250)(-138,250)(-150,265)
\bezier50(-178,300)(-178,320)(-175,325)
\bezier50(-178,300)(-178,280)(-175,275)
\bezier5(-172,300)(-172,320)(-175,325)
\bezier5(-172,300)(-172,280)(-175,275)
\bezier100(-131,300)(-131,340)(-125,350)
\bezier100(-131,300)(-131,260)(-125,250)
\bezier10(-119,300)(-119,340)(-125,350)
\bezier10(-119,300)(-119,260)(-125,250)
\bezier70(-154,300)(-154,328)(-150,335)
\bezier70(-154,300)(-154,272)(-150,265)
\bezier7(-146,300)(-146,328)(-150,335)
\bezier7(-146,300)(-146,272)(-150,265)
\put(-302,300){\line(1,0){22}}
\put(-176,300){\line(1,0){20}}
\put(-152,300){\line(1,0){19}}
\put(-129,300){\line(1,0){23}}
\put(-102,300){\line(1,0){22}}
\bezier100(-75,325)(-88,325)(-100,335)
\bezier100(-125,350)(-112,350)(-100,335)
\bezier100(-75,275)(-88,275)(-100,265)
\bezier100(-125,250)(-112,250)(-100,265)
\bezier5(-72,300)(-72,320)(-75,325)
\bezier5(-72,300)(-72,280)(-75,275)
\bezier50(-78,300)(-78,320)(-75,325)
\bezier50(-78,300)(-78,280)(-75,275)
\bezier7(-96,300)(-96,328)(-100,335)
\bezier7(-96,300)(-96,272)(-100,265)
\bezier70(-104,300)(-104,328)(-100,335)
\bezier70(-104,300)(-104,272)(-100,265)
\bezier100(-275,325)(-262,325)(-250,335)
\bezier100(-225,350)(-238,350)(-250,335)
\bezier100(-275,275)(-262,275)(-250,265)
\bezier100(-225,250)(-238,250)(-250,265)
\bezier50(-278,300)(-278,320)(-275,325)
\bezier50(-278,300)(-278,280)(-275,275)
\bezier5(-272,300)(-272,320)(-275,325)
\bezier5(-272,300)(-272,280)(-275,275)
\bezier100(-231,300)(-231,340)(-225,350)
\bezier100(-231,300)(-231,260)(-225,250)
\bezier10(-219,300)(-219,340)(-225,350)
\bezier10(-219,300)(-219,260)(-225,250)
\bezier70(-254,300)(-254,328)(-250,335)
\bezier70(-254,300)(-254,272)(-250,265)
\bezier7(-246,300)(-246,328)(-250,335)
\bezier7(-246,300)(-246,272)(-250,265)
\put(-276,300){\line(1,0){20}}
\put(-252,300){\line(1,0){19}}
\put(-229,300){\line(1,0){23}}
\put(-202,300){\line(1,0){22}}
\bezier100(-175,325)(-188,325)(-200,335)
\bezier100(-225,350)(-212,350)(-200,335)
\bezier100(-175,275)(-188,275)(-200,265)
\bezier100(-225,250)(-212,250)(-200,265)
\bezier7(-196,300)(-196,328)(-200,335)
\bezier7(-196,300)(-196,272)(-200,265)
\bezier70(-204,300)(-204,328)(-200,335)
\bezier70(-204,300)(-204,272)(-200,265)
\bezier100(-75,325)(-62,325)(-50,335)
\bezier100(-25,350)(-38,350)(-50,335)
\bezier100(-75,275)(-62,275)(-50,265)
\bezier100(-25,250)(-38,250)(-50,265)
\bezier100(-31,300)(-31,340)(-25,350)
\bezier100(-31,300)(-31,260)(-25,250)
\bezier10(-19,300)(-19,340)(-25,350)
\bezier10(-19,300)(-19,260)(-25,250)
\bezier70(-54,300)(-54,328)(-50,335)
\bezier70(-54,300)(-54,272)(-50,265)
\bezier7(-46,300)(-46,328)(-50,335)
\bezier7(-46,300)(-46,272)(-50,265)
\put(-76,300){\line(1,0){20}}
\put(-52,300){\line(1,0){19}}
\put(-29,300){\line(1,0){23}}
\put(-2,300){\line(1,0){22}}
\bezier100(25,325)(12,325)(0,335)
\bezier100(-25,350)(-12,350)(0,335)
\bezier100(25,275)(12,275)(0,265)
\bezier100(-25,250)(-12,250)(0,265)
\bezier5(28,300)(28,320)(25,325)
\bezier5(28,300)(28,280)(25,275)
\bezier50(22,300)(22,320)(25,325)
\bezier50(22,300)(22,280)(25,275)
\bezier7(4,300)(4,328)(0,335)
\bezier7(4,300)(4,272)(0,265)
\bezier70(-4,300)(-4,328)(0,335)
\bezier70(-4,300)(-4,272)(0,265)
\bezier100(25,325)(38,325)(50,335)
\bezier100(75,350)(62,350)(50,335)
\bezier100(25,275)(38,275)(50,265)
\bezier100(75,250)(62,250)(50,265)
\bezier100(69,300)(69,340)(75,350)
\bezier100(69,300)(69,260)(75,250)
\bezier100(81,300)(81,340)(75,350)
\bezier100(81,300)(81,260)(75,250)
\bezier70(46,300)(46,328)(50,335)
\bezier70(46,300)(46,272)(50,265)
\bezier7(54,300)(54,328)(50,335)
\bezier7(54,300)(54,272)(50,265)
\put(24,300){\line(1,0){20}}
\put(48,300){\line(1,0){19}}
\put(71,300){\line(1,0){45}}
\end{picture}
\vspace{0.75in}
\end{center}
    \caption{A portion of ${\cal D}(\mu) \subset \R^3$, $\mu > 0$.}%
    \label{fig:delaunayR3}
\end{figure}


\paragraph{3.2 Special parts of Euclidean Delaunay
unduloids $\cD(\mu)$.}\label{SS:spdh}
Without loss of generality,
we may assume that the function $f$ defining the profile curve of
$\cD(\mu)$ satisfies $f(0) = a_-(\mu)$. It follows easily that $f(T(\mu)) =
a_-(\mu), f(T(\mu)/2) = a_+(\mu)$ and that $f$ is symmetric with
respect to the values $k\, T(\mu)/2, \, k \in \Z$.

\smallskip

Let the \emph{basic Dirichlet block} for $\cD(\mu)$ be the compact domain

\begin{equation}\label{E:bbdh-1}
\cB(\mu) := F([0, \frac{T(\mu)}{2}]\times S^{n-1}) \mbox{~~or~~}
F([\frac{T(\mu)}{2}, T(\mu)]\times S^{n-1})\,,
\end{equation}

where $F$ is the parametrization (\ref{E:rs1}), see Figure 1.
We also introduce the pieces $\cB_{\ell}(\mu)$ obtained by glueing
$\ell$ basic Dirichlet blocks,

\begin{equation}\label{E:bbdh-1a}
\cB_{\ell}(\mu) := F([0, \ell \,\frac{T(\mu)}{2}]\times S^{n-1})
\mbox{~~or~~} F([\frac{T(\mu)}{2}, (\ell+1) \,
\frac{T(\mu)}{2}]\times S^{n-1})\; .
\end{equation}

\smallskip

Let $a$ be the function

\begin{equation}\label{E:bbdh-2}
a (x) = \ip{N(x,\omega)}{(1,0,...,0)}= (1+f'^{2}(x))^{-1/2}\, f'(x) \;
,
\end{equation}

where $f$ is the profile curve of the Euclidean Delaunay unduloid $\cD(\mu)$.

\begin{lemma}\label{L:bbe}
The function $a$ satisfies $(\Delta - V) a = 0$ and vanishes
precisely at the half-integer multiples of $T(\mu)$. Furthermore,
$a'$ has exactly two zeroes $\zeta_1(\mu), \zeta_2(\mu)$ in the
interval $[0, T(\mu)]$, with

$$ 0 < \zeta_1(\mu) < \frac{T(\mu)}{2} < \zeta_2(\mu) < T(\mu) \; .$$

(For a true cylinder, $a(x)$ is constant, so the
values $\zeta_1(\frac{1}{n}(\frac{n-1}{n})^{n-1})$ and
$\zeta_2(\frac{1}{n}(\frac{n-1}{n})^{n-1})$ must be determined by
the limits of $\zeta_1(\mu)$ and $\zeta_2(\mu)$ as $\mu$ increases
to $\frac{1}{n}(\frac{n-1}{n})^{n-1}$.)
\end{lemma}

\pf The first assertion is classical: the scalar product of the
unit normal vector of a cmc hypersurface with a Killing field is a
solution of $L u = 0$ (see \cite{choe}, page 196, or the proof
of Theorem 2.7 in \cite{bgs}). The second assertion is obvious. The
assertion on the zeroes of $a'$ follows from the fact that
$(\Delta - V)a=0$ reduces to a Sturm--Liouville equation, since
the functions $a$ and $V$ depend on the variable $x$ only. \qedd

\medskip

Let the \emph{basic Neumann block} for $\cD(\mu)$ be the compact domain

\begin{equation}\label{E:bbdh-3}
\cC(\mu) := F([\zeta_1(\mu), T(\mu) + \zeta_1(\mu)]\times
S^{n-1})\; .
\end{equation}

\smallskip

We also introduce the pieces $\cC_{\ell}(\mu)$ obtained by glueing
$\ell$ basic Neumann blocks, see Figure~1,

\begin{equation}\label{E:bbdh-3a}
\cC_{\ell}(\mu) := F([\zeta_1(\mu), \ell \, T(\mu) +
\zeta_1(\mu)]\times S^{n-1})\; .
\end{equation}


\bigskip


\begin{figure}[h] 
\begin{center}
\unitlength=0.7pt
\begin{picture}(-200.00,200.00)(0.00,0.00)
\put(-200.00,25.00){\makebox(0,0)[cc]{\small profile curve}}
\put(-135,12){\makebox(0,0)[cc]{\small $\varphi(t)$}}
\put(-28,35){\makebox(0,0)[cc]{\small $m(t)$}}
\put(-58,70){\makebox(0,0)[cc]{\small $\rho(t)$}}
\put(-93,106){\makebox(0,0)[cc]{\small $\gamma$}}
\put(-110,76){\makebox(0,0)[cc]{\small $e^t$}}
\put(-100,-8){\makebox(0,0)[cc]{\small $0$}}
\put(-172,28){\vector(1,1){11}} \put(-123,12){\vector(1,0){30}}
\put(-100,0){\line(0,1){170}} \put(-100,0){\line(1,1){160}}
\put(-100,0){\line(1,0){160}} \put(-100,0){\line(-1,0){160}}
\put(-100,0){\line(-1,1){160}} \bezier20(-100,71)(-54,68)(-40,40)
\bezier20(-100,7.1)(-95.4,6.8)(-94,4)
\bezier30(-100,0)(-70,20)(-40,40)
\bezier15(-85,10)(-82.5,12.5)(-82.5,17.5)
\bezier40(-80,30)(-82.5,27.5)(-82.5,17.5)
\bezier15(-115,10)(-117.5,12.5)(-117.5,17.5)
\bezier40(-120,30)(-117.5,27.5)(-117.5,17.5)
\bezier30(-80,30)(-75,35)(-65,35)
\bezier80(-40,40)(-45,35)(-65,35)
\bezier60(-40,40)(-30,50)(-30,70)
\bezier160(-20,120)(-30,110)(-30,70)
\bezier120(-20,120)(0,140)(40,140)
\bezier30(-120,30)(-125,35)(-135,35)
\bezier80(-160,40)(-155,35)(-135,35)
\bezier60(-160,40)(-170,50)(-170,70)
\bezier160(-180,120)(-170,110)(-170,70)
\bezier120(-180,120)(-200,140)(-240,140)
\bezier120(-100,125)(25,125)(40,140)
\bezier120(-100,125)(-225,125)(-240,140)
\bezier20(-100,155)(25,155)(40,140)
\bezier20(-100,155)(-255,155)(-240,140)
\end{picture}
\vspace{0.2in}
\end{center}
    \caption{A portion of ${\cal D}_H(\mu) \subset \HH^3$, $\mu > 0$,
             $H > 1$.}%
    \label{fig:delaunayH3}
\end{figure}
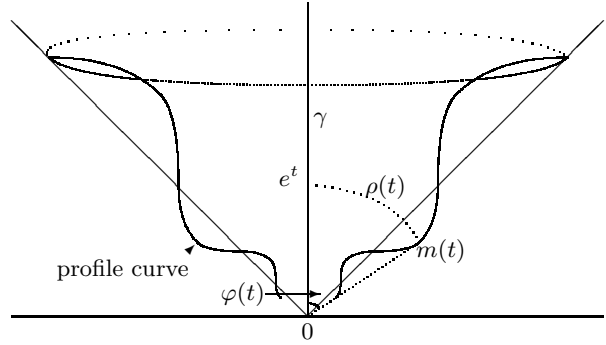


\paragraph{3.3 Delaunay unduloids with
cmc $H>1$ in hyperbolic space.}\label{SS:dhhs} We choose the
half-space model $\big\aco (x_1, \ldots , x_n,y) \in \R^{n+1}
~\big|~ y > 0\big\acf$ for hyperbolic space $\HH^{n+1}$ (with the
hyperbolic space metric), and we fix the geodesic $\gamma(t)= (0,
\ldots , 0, e^{t})$.

\smallskip

The profile curve of a hyperbolic Delaunay unduloid is described,
say in the vertical $2$-dimensional plane $\big\aco
x_1,y\big\acf$, as a geodesic graph. The point $m(t)$ on the
profile curve  is at geodesic distance $\rho(t)$ from the point
$\gamma(t)$. Let $\varphi(t)$ be the angle $\angle (\gamma(t)\, 0
\, m(t))$, see Figure~2. Then, $\sinh \rho(t) = \tan \varphi(t)$.

\smallskip

With these notations, the profile curve is given by $\big( e^{t}
\sin \varphi(t), e^{t} \cos \varphi(t)\big)$, where $\varphi$
satisfies the differential equation (\cite{KKMS}, Equation (6.3)
page 34)

\begin{equation}\label{E:del-h2}
\mu =
\frac{(\tan \varphi)^{n-1}}{\cos \varphi \sqrt{1+\varphi'^{2}}} -
H (\tan \varphi)^{n}\,.
\end{equation}

Here, $\mu > 0$ is the weight parameter, and the
(normalized) mean curvature $H$ satisfies $H > 1$.  (Note that the
mean curvature is not normalized in \cite{KKMS}.)
The hyperbolic Delaunay unduloids $\cD_H(\mu)$ are given by

\begin{equation}\label{E:del-h1}
\R \times S^{n-1} \ni (t,\omega) \stackrel{\Phi}{\to} (e^{t} \sin
\varphi(t) \, \omega, e^{t} \cos \varphi(t)) =: (f(t) \, \omega,
g(t))\in \HH^{\,n+1}\; .
\end{equation}
As in the case of the Euclidean Delaunay unduloids, it can be
shown that the function $\varphi$ (or equivalently $\rho$) is
pinched between two values $0 < \alpha_-(\mu) \le \varphi(t) \le
\alpha_+(\mu)$ and periodic with period $\tau(\mu)$. The Delaunay
unduloids obtained in this way with $\mu > 0$ are embedded.

\smallskip

A unit normal vector to the hypersurface $\cD_H(\mu)$ is given (with
the above notations) by

\begin{equation}\label{E:del-h3}
N(t,\omega) = \frac{\cos \varphi}{\sqrt{1 + \varphi'^{2}}} \, (g'
\, \omega, - f')\,.
\end{equation}

\medskip

The stability operator $L$ is of the form

\begin{equation}\label{E:del-h5}
        L = \Delta - V \; , \;\;\;\;\;\;\;\;\;\; \mbox{where} \;\;
        V = -n + \|B\|^2 \; .
\end{equation}

Note that $V$ is periodic and hence bounded on $\R$, as in the
Euclidean case. There is a nice expression for the function $V$ in
the hyperbolic case:

\begin{lemma}\label{L:new2}
  With the notations as in Equations (\ref{E:del-h1}) and
(\ref{E:del-h5}), we have, for the $n$-dimensional hyperbolic Delaunay
unduloid $\cD_H(\mu)$, that
$$
V(\Phi(t,\omega)) = n(H^2 - 1) + n (n-1) \mu^2 (\tan
\varphi)^{-2n} \; .
$$
\end{lemma}

\pf For $n=2$ a proof is in \cite{C}.
The case $n\ge 3$ is similar, using computations
in \cite{hsiang}. \qedd

\medskip

In order to estimate the index of certain pieces of $\cD_H(\mu)$,
we need the following lemma.

\begin{lemma}\label{L:new3}
With the notations as in Equations (\ref{E:del-h5}) and
(\ref{E:del-h2}), there exists a constant $c(n,H)$ depending only on
$n$ and $H$ so that
$$
V \tan^2 \varphi \le c(n,H)
$$
on $\cD_H(\mu)$ with weight parameter $\mu$ and cmc $H > 1$.
\end{lemma}

\pf The proof is left to the reader (use Equation
(\ref{E:del-h2}), Lemma \ref{L:new2}, and compute).\qedd

\paragraph{3.4 Special parts of hyperbolic Delaunay unduloids $\cD_H(\mu)$.}
\label{SS:sphdh}
Without loss of
generality, we may assume $\varphi$ satisfies $\varphi(0) =
\alpha_-(\mu)$. Thus $\varphi(\tau(\mu)) = \alpha_-(\mu),
\varphi(\tau(\mu)/2) = \alpha_+(\mu)$ and $\varphi$ is
symmetric with respect to the values $k\, \tau(\mu)/2, \, k \in
\Z$.

\smallskip

Analogous to the Euclidean case, we define the
\emph{basic Dirichlet block} $\cB(\mu)$ for
$\cD_H(\mu)$, the glueing of $\ell$ basic Dirichlet blocks
$\cB_{\ell}(\mu)$, and the function $a$ as:
\begin{equation}\label{E:bbhdh-1}
\cB(\mu) := \Phi([0, \frac{\tau(\mu)}{2}]\times S^{n-1})
\mbox{~~or~~} \Phi([\frac{\tau(\mu)}{2}, \tau(\mu)]\times S^{n-1})
\; ,
\end{equation}
\begin{equation}\label{E:bbhdh-1a}
\cB_{\ell}(\mu) := \Phi([0, \ell \,\frac{\tau(\mu)}{2}]\times
S^{n-1}) \mbox{~~or~~} \Phi([\frac{\tau(\mu)}{2},
(\ell+1) \, \frac{\tau(\mu)}{2}]\times S^{n-1})\; ,
\end{equation}
\begin{equation}\label{E:bbhdh-2}
a (x) = \ip{N(t,\omega)}{\mathcal{Y}} =\frac{\varphi'(x)}{\cos
\varphi(x) \, \sqrt{1 + \varphi'^{2}(x)}} \; ,
\end{equation}
where $\Phi$ is the parametrization (\ref{E:del-h1}), and
$\varphi$ satisfies (\ref{E:del-h2}), and $\mathcal{Y}$ is the
Killing field corresponding to hyperbolic translation along
the axis of the Delaunay unduloid.  The following lemma is proved in
the same way as Lemma \ref{L:bbe}:

\begin{lemma}\label{L:bbh}
The function $a$ satisfies $(\Delta - V) a = 0$ and
vanishes precisely at the half-integer multiples of $\tau(\mu)$.
Furthermore, $a'$ has exactly two zeroes $\zeta_1(\mu),
\zeta_2(\mu)$ in the interval $[0, \tau(\mu)]$, with
$$ 0 < \zeta_1(\mu) < \frac{\tau(\mu)}{2} < \zeta_2(\mu) < \tau
(\mu).$$

(Again, the values $\zeta_j(\mu)$ and $\tau(\mu)$ for the true
hyperbolic cylinder are determined as limiting values of the
$\zeta_j(\mu)$ and $\tau(\mu)$ for noncylindrical hyperbolic
Delaunay unduloids.)
\end{lemma}

Let the \emph{basic Neumann block} $\cC(\mu)$ for
$\cD_H(\mu)$ and the glueing of
$\ell$ basic Neumann blocks $\cC_{\ell}(\mu)$ be

\begin{equation}\label{E:bbhdh-3}
\cC(\mu) := \Phi([\zeta_1(\mu), \tau(\mu) + \zeta_1(\mu)]\times
S^{n-1})\; ,
\end{equation}
\begin{equation}\label{E:bbhdh-3a}
\cC_{\ell}(\mu) := \Phi([\zeta_1(\mu), \ell \, \tau(\mu) +
\zeta_1(\mu)]\times S^{n-1})\; .
\end{equation}

\section{Index estimates for pieces of Delaunay unduloids}\label{S:iecpdh}

\paragraph{4.1 Preliminary results.}\label{Par:pr} We state the following lemma
for later reference.

\medskip

\begin{lemma}\label{L:l5}
Let $A, B, V$ be smooth bounded functions on $\R$. Assume that $A, B$ are
bounded from below by a positive constant. Let $P$ be the manifold
$[a,b]\times S^{n-1}$ equipped
with the metric $g := A^2(x) dx^2 + B^2(x) g_S$, where $g_S$ is
the canonical metric on $S^{n-1}$. We are interested in the
eigenvalue problem $(\Delta_g - V)\, y(x,\omega) = \lambda\,
y(x,\omega)$, with Dirichlet or Neumann conditions on
$\partial P = (\aco a \acf \times S^{n-1}) \cup
(\aco b \acf \times S^{n-1})$. Let $\Lambda_k = k(k+n-2), \,
k \ge 0$, denote the eigenvalues of the Laplacian on $S^{n-1}$ and
let $m(\Lambda_k)$ denote the multiplicity of $\Lambda_k$ (this is
a polynomial in $k$, of degree $n-2$). Let $L := \Delta_g - V$
and define the operators $L_k, k\ge 0,$ by
$$ L_k \, u= - \frac{d}{dx}\left( A^{-1}B^{n-1} \,
\frac{du}{dx}\right) + A B^{n-3} \, \left( \Lambda_k - B^{2}\,
V\right)\, u\; . $$ Let us denote by $\sigma (L)$ the set of
eigenvalues of $L$, counted with multiplicities, and by
$\sigma(L_k)$ the eigenvalues of the problem $L_k \, u = \lambda
\, AB^{n-1}\, u$. Then
$$ \sigma(L) = \bigsqcup_{k=0}^{\infty} m(\Lambda_k)
\,\sigma(L_k) \; , $$
where the expression in the right-hand side means that each
eigenvalue of $L_k$ appears with multiplicity $m(\Lambda_k)$ in
$\sigma(L)$ (summing up multiplicities if the same number
$\lambda$ appears in several $\sigma(L_k)$). In particular, the
index (number of negative eigenvalues) of $L$ is given by
$$ \mathrm{Index}(L) = \sum_{k=0}^{\infty}\mathrm{Index}(L_k) \; , $$
and the sum on the right-hand side involves only finitely many terms.
\end{lemma}

\pf If $y(x,\omega)$ satisfies $L y = \lambda y$, then
$$ - \frac{\partial}{\partial x}\left( A^{-1}B^{n-1} \,
\frac{\partial y}{\partial x}\right) + A B^{n-3}\, \left( \Delta_S
\,y - B^{2}\, V y\right) = \lambda \, AB^{n-1} y\,.$$

In order to prove the lemma, it suffices to decompose the function
$y(x,\omega)$ into a series of spherical harmonics. The generic
term in this series will be of the form $u(x) Y(\omega)$ where $Y$
is a $k$-spherical harmonic and the preceding equation becomes
$$ - \frac{d}{dx}\left( A^{-1}B^{n-1} \, \frac{du}{dx}\right) + A
B^{n-3}\, \left(\Lambda_k - B^{2}\, V\right)\,u = \lambda \,
AB^{n-1} u\,.$$

The first assertion of the lemma follows easily. For the final
assertion we need only remark that $\Lambda_k$ tends to
infinity with $k$ and hence the operators $L_k$ are positive
for $k$ large enough (this is because $A, B$ are bounded from
below by positive constants and $V$ is bounded). \qedd

\medskip

\textbf{Remark.~} The $L_k$ operate on functions of a single
variable, hence their eigenvalues are all simple, and an
eigenfunction associated to the $j$'th eigenvalue has exactly $j$
nodal domains.  We use these properties in upcoming arguments.

\medskip

Let $\cD$ be a Delaunay unduloid in $\R^{\, n+1}$ (with $H = 1,
\mu > 0$) or  in $\HH^{\, n+1}$ (with $H > 1, \mu > 0$), with cmc $H$
and weight parameter $\mu$. Let $\R^{\,
n+1}_+$ and $\HH^{\, n+1}_+$ denote one of the closed half-spaces
defined by a geodesic hyperplane containing the axis of $\cD$.
Then we have:

\begin{proposition}\label{P:sshdp}
The stability operator $\Delta - V$ of the Delaunay unduloid $\cD$
is positive in any $\Omega$ contained in $\cD_+ := \cD \bigcap
\R^{\, n+1}_+$ or $\cD_+ := \cD \bigcap \HH^{\,n+1}_+$, with
respect to Dirichlet boundary conditions. In particular, the
half-Delaunay unduloids $\cD_+$ are (strongly) stable.
\end{proposition}

\pf This result is well-known for Euclidean graphs.

\smallskip

In the hyperbolic case, one has to be more careful, as certain kinds
of graphs are not stable.  To prove the proposition, it
suffices to find a positive solution of $(\Delta - V)
\, y = 0 $ on $\cD_+$. Such a solution will be given by the normal
component of a well chosen Killing field.

We consider the Killing field $\mathcal{Y}_{\theta}(\omega,t) =
(\theta,0)$ in $\HH^{\; n+1 }$, where $\theta \in S^{n-1}$ is chosen so
that $\mathcal{Y}_{\theta}(\omega,t)|_{\partial \cD_+}$ is perpendicular to
the geodesic hyperplane containing $\partial \cD_+$.
The function
$$a_{\theta}(t,\omega) := \ip{N}{\mathcal{Y}_{\theta}}$$
satisfies $(\Delta - V)a_\theta = 0$ and
is equal to $g'(t)\, \ip{\theta}{\omega}$ up to a positive factor
(recall that $g(t) = e^t \cos \varphi(t)$, see Equation (\ref{E:del-h1})).
To prove that $a_{\theta}>0$ in the interior of $\cD_+$, it suffices to look
at the sign of $g'(t) =
e^{t}\big( \cos \varphi(t) - \varphi'(t) \, \sin \varphi(t)\big)$.
Assume there is a point $t_0$ at which $g'$ vanishes, then
\[ {1 \over \varphi^\prime(t_0)} = \tan \varphi(t_0) > 0 \]
and Equation (\ref{E:del-h2}) implies
$0 < \mu = (1-H)/((\varphi^\prime(t_0))^n) < 0$, a contradiction.
\qedd

\paragraph{4.2 Index estimates for certain Delaunay
pieces.}\label{SS:iecdp}
We have the following estimates for the indexes of
the Delaunay pieces $\cB_{\ell}$ and $\cC_{\ell}$, in both Euclidean
and hyperbolic cases:

\begin{proposition}\label{P:ibp}
The index of the Delaunay piece $\cB_{\ell}(\mu)$ with Dirichlet
conditions at both boundary components is exactly $\ell - 1$.
\end{proposition}

\begin{proposition}\label{P:inp}
There is a constant $c_1(n,H)$, which depends only on the
dimension $n$ and the mean curvature $H$, such that the index of
the Delaunay piece $\cC_{\ell}(\mu)$ with Neumann conditions at
both boundary components satisfies

$$ 2 \ell \le \mbox{Neumann~Index}(\cC_{\ell}(\mu)) \le 2 \ell
+ c_1(n,H)\,.$$
\end{proposition}

\textbf{Proofs.~} The proofs of these two propositions are quite
similar.

\textbf{Step 1.~} The induced metric on the pieces $\cB_{\ell}(\mu)$ or
$\cC_{\ell}(\mu)$ is of
the type described in Lemma \ref{L:l5}, with $A=\sqrt{1+(f^\prime)^2}$,
$B=f$ in the Euclidean case, and $A=\sqrt{(1+(\varphi^\prime)^2)(1+
\tan^2 \varphi)}$, $B=\tan \varphi$ in the hyperbolic case.
Hence, to estimate the index we only need to look at the
indexes of the corresponding operators $A^{-1} B^{1-n} L_k$,
and we already know that for $k$ large enough the
operator $A^{-1} B^{1-n} L_k$ is positive, implying that its index is zero.
In fact, looking at the bounds we have for $B^2 V$ in Lemmas \ref{L:new1} and
\ref{L:new3}, we see that there exists a constant $c(n,H)$ such that
$A^{-1} B^{1-n} L_k$ is positive whenever $k\ge c(n,H)$.

\medskip

\textbf{Proof of Proposition \ref{P:ibp}, Step 2.~} The following
proof applies when the Delaunay unduloid is not a cylinder. (For
cylinders the index estimates are trivial, and we do not include
the arguments here.) By Lemmas \ref{L:bbe} and \ref{L:bbh}, the
function $a$ in (\ref{E:bbdh-2}) and (\ref{E:bbhdh-2}) satisfies
$(\Delta - V) a = 0$ and $a\left|_{\partial \cB_{\ell}} = 0
\right.$. This function $a$ has precisely $\ell$ nodal domains in
$\cB_{\ell}$. It follows that $0$ is the $\ell$-th eigenvalue of
the operator $A^{-1} B^{1-n} L_0$. So the Dirichlet index of
$\cB_{\ell}$ is at least $\ell - 1$ and is bigger than $\ell - 1$
if and only if some of the operators $A^{-1} B^{1-n} L_k, k\ge 1$,
have negative eigenvalues. Assume this is the case and that some
$u(x)$ with Dirichlet boundary conditions satisfies $A^{-1}
B^{1-n} L_k \, u = \lambda \, u$ for some $\lambda < 0$. This
implies that

$$ (\Delta - V) \, u \, Y = \lambda \, u \, Y$$

for any spherical harmonic $Y$ of degree $k$. Choosing, for example, a
radial spherical harmonic $Y$, we can always find a domain $D_k
\subset \cD_+$ such that the function $u \, Y$ is positive in
$D_k$ and satisfies

$$
\left\aco
\begin{array}{ll}
  (\Delta - V) \, u \, Y = \lambda \, u \, Y & \mbox{in~} D_k\,, \\
  u\, Y = 0 & \mbox{on~} \partial D_k\,.
\end{array}
\right.
$$

This contradicts Proposition \ref{P:sshdp}.

\medskip

\textbf{Proof of Proposition \ref{P:inp}, Step 2.~} We use the
same function $a$ as before, so $(\Delta - V) a = 0$. The domain
$\cC_{\ell}$ was designed so that $a' = \frac{\partial a}{\partial
\nu} = 0$ on $\partial \cC_{\ell}$ and so that $a$ has exactly
$(2\ell +1)$ nodal domains in $\cC_{\ell}$. It follows, as in the
preceding argument, that $0$ is the $(2\ell + 1)$-st eigenvalue of
$L_0$ and hence that the Neumann index of $\cC_{{\ell}}$ is at
least $2\ell$. In order to obtain the upper bound, we remark that
the Neumann index of $A^{-1} B^{1-n} L_k, k\ge 1,$ is at most $2$.
Indeed, assume it is at least $3$. Then there is an eigenfunction
$u$ of $A^{-1} B^{1-n} L_k, k\ge 1,$ with at least three nodal
domains and hence with an interior nodal domain. We can then
repeat the argument in Step 2, proof of Proposition \ref{P:ibp},
and arrive at a contradiction to Proposition \ref{P:sshdp}. An
eigenvalue of $A^{-1} B^{1-n} L_k$ gives an eigenvalue of $L$ with
multiplicity a polynomial of degree $(n-2)$ in $k$. Since $A^{-1}
B^{1-n} L_k$ is positive for $k \ge c(n,H)$, the result follows.
\qedd

\section{Index growth results}\label{S:igr}

\paragraph{5.1 Eigenvalue estimates for almost Delaunay
pieces.}\label{SS:eeadp}
Fix a Delaunay unduloid $\cD$ and a piece $\cE \subset \cD$
which is bounded by two ``parallel spheres" in geodesic
hyperplanes orthogonal to the axis of revolution. We call
$\widetilde{\cE}$ an \emph{almost Delaunay piece} if it is a
cylindrical graph over $\cE$.

\begin{lemma}\label{L:l10}
There exists a constant $c_2(n,H)$, depending only on the
dimension $n$ and mean curvature $H$, such that if
$\widetilde{\cE}$ is close enough to $\cE$ in the $C^{2}$-sense,
then

$$ \ind(\cE) \le \ind(\widetilde{\cE}) \le \ind(\cE) +
c_2(n,H),$$

where $\ind$ denotes the index for either Dirichlet or Neumann
conditions on the corresponding boundary components of $\partial
\cE, \partial \widetilde{\cE}.$

\end{lemma}

\pf Indeed, once the piece $\cE$ is \emph{fixed}, we can write the
eigenvalues of the operator $L$ on $\cE$ (with respect to some
Dirichlet or Neumann conditions on the boundary components) as

$$ \lambda_1(\cE) < \lambda_2(\cE) \le \ldots \le \lambda_k(\cE) <
0 \le \lambda_{k+1}(\cE) \le \ldots $$

where $k = \ind(\cE)$. If $\widetilde{\cE}$ is close enough to
$\cE$ in the $C^{2}$-sense, the negative eigenvalues of the
operator $\widetilde{L}$ corresponding to $L$ are close to the
corresponding eigenvalues of $L$. It follows that $\ind(\cE) \le
\ind(\widetilde{\cE})$ because $\lambda_k(\widetilde{\cE}) < 0$,
with $\ind(\cE) = \ind(\widetilde{\cE})$ unless
$\lambda_{k+1}(\cE) = 0$, in which case we may have
$\lambda_{k+1}(\widetilde{\cE}) < 0$ and the constant $c_2(n,H)$
takes the possible multiplicity of $\lambda_{k+1}(\cE)$ into
account. This multiplicity can be bounded as indicated in the
proof of Proposition \ref{P:inp}. \qedd

\medskip \textbf{Note.~} On $\cE$, the eigenvalue problem $L_p u =
\lambda A B^{n-1} u$ introduced earlier
cannot have eigenvalue $0$ when $p\ge 1$ in case of Dirichlet
boundary condition (by arguing like in the proof of
Proposition \ref{P:ibp}).


\bigskip

\begin{figure}[h] 
\begin{center}
\unitlength=0.7pt
\begin{picture}(-200.00,200.00)(0.00,200.00)
\put(-114.00,200.00){\makebox(0,0)[cc]{\small $M_0$}}
\put(-175,165){\makebox(0,0)[cc]{\small $E_1^R$}}
\put(-10,155){\makebox(0,0)[cc]{\small $E_2^R$}}
\put(-10,245){\makebox(0,0)[cc]{\small $E_3^R$}}
\put(-130,200){\vector(-1,0){140}} \put(-130,190){\line(0,1){20}}
\bezier20(-130,210)(-135,210)(-140,215)
\bezier20(-150,220)(-145,220)(-140,215)
\bezier20(-170,210)(-165,210)(-160,215)
\bezier20(-150,220)(-155,220)(-160,215)
\bezier20(-170,210)(-175,210)(-180,215)
\bezier20(-190,220)(-185,220)(-180,215)
\bezier20(-210,210)(-205,210)(-200,215)
\bezier20(-190,220)(-195,220)(-200,215)
\bezier20(-210,210)(-215,210)(-220,215)
\bezier20(-230,220)(-225,220)(-220,215)
\bezier20(-250,210)(-245,210)(-240,215)
\bezier20(-230,220)(-235,220)(-240,215)
\bezier20(-130,190)(-135,190)(-140,185)
\bezier20(-150,180)(-145,180)(-140,185)
\bezier20(-170,190)(-165,190)(-160,185)
\bezier20(-150,180)(-155,180)(-160,185)
\bezier20(-170,190)(-175,190)(-180,185)
\bezier20(-190,180)(-185,180)(-180,185)
\bezier20(-210,190)(-205,190)(-200,185)
\bezier20(-190,180)(-195,180)(-200,185)
\bezier20(-210,190)(-215,190)(-220,185)
\bezier20(-230,180)(-225,180)(-220,185)
\bezier20(-250,190)(-245,190)(-240,185)
\bezier20(-230,180)(-235,180)(-240,185)
\put(-75,225){\vector(1,1){70}} \put(-70,220){\line(-1,1){10}}
\bezier15(-60,250)(-62,248)(-67,247)
\bezier15(-75,245)(-73,247)(-67,247)
\bezier15(-50,240)(-52,238)(-53,233)
\bezier15(-55,225)(-53,227)(-53,233)
\bezier15(-80,230)(-78,232)(-77,237)
\bezier15(-75,245)(-77,243)(-77,237)
\bezier15(-70,220)(-68,222)(-63,223)
\bezier15(-55,225)(-57,223)(-63,223)
\bezier15(-40,270)(-42,268)(-47,267)
\bezier15(-55,265)(-53,267)(-47,267)
\bezier15(-30,260)(-32,258)(-33,253)
\bezier15(-35,245)(-33,247)(-33,253)
\bezier15(-60,250)(-58,252)(-57,257)
\bezier15(-55,265)(-57,263)(-57,257)
\bezier15(-50,240)(-48,242)(-43,243)
\bezier15(-35,245)(-37,243)(-43,243)
\bezier15(-20,290)(-22,288)(-27,287)
\bezier15(-35,285)(-33,287)(-27,287)
\bezier15(-10,280)(-12,278)(-13,273)
\bezier15(-15,265)(-13,267)(-13,273)
\bezier15(-40,270)(-38,272)(-37,277)
\bezier15(-35,285)(-37,283)(-37,277)
\bezier15(-30,260)(-28,262)(-23,263)
\bezier15(-15,265)(-17,263)(-23,263)
\put(-75,175){\vector(1,-1){70}} \put(-70,180){\line(-1,-1){10}}
\bezier15(-60,150)(-62,152)(-67,153)
\bezier15(-75,155)(-73,153)(-67,153)
\bezier15(-50,160)(-52,162)(-53,167)
\bezier15(-55,175)(-53,173)(-53,167)
\bezier15(-80,170)(-78,168)(-77,163)
\bezier15(-75,155)(-77,157)(-77,163)
\bezier15(-70,180)(-68,178)(-63,177)
\bezier15(-55,175)(-57,177)(-63,177)
\bezier15(-40,130)(-42,132)(-47,133)
\bezier15(-55,135)(-53,133)(-47,133)
\bezier15(-30,140)(-32,142)(-33,147)
\bezier15(-35,155)(-33,153)(-33,147)
\bezier15(-60,150)(-58,148)(-57,143)
\bezier15(-55,135)(-57,137)(-57,143)
\bezier15(-50,160)(-48,158)(-43,157)
\bezier15(-35,155)(-37,157)(-43,157)
\bezier15(-20,110)(-22,112)(-27,113)
\bezier15(-35,115)(-33,113)(-27,113)
\bezier15(-10,120)(-12,122)(-13,127)
\bezier15(-15,135)(-13,133)(-13,127)
\bezier15(-40,130)(-38,128)(-37,123)
\bezier15(-35,115)(-37,117)(-37,123)
\bezier15(-30,140)(-28,138)(-23,137)
\bezier15(-15,135)(-17,137)(-23,137)
\bezier80(-70,220)(-90,200)(-70,180)
\bezier100(-130,210)(-100,210)(-80,230)
\bezier100(-130,190)(-100,190)(-80,170)
\bezier50(-65,185)(-62.5,187.5)(-60,190)
\bezier50(-5,125)(-2.5,127.5)(0,130)
\bezier50(-65,215)(-62.5,212.5)(-60,210)
\bezier50(-5,275)(-2.5,272.5)(0,270)
\put(-60,190){\line(1,-1){60}} \put(-60,210){\line(1,1){60}}
\put(-130,175){\line(0,1){5}} \put(-220,175){\line(0,1){5}}
\put(-130,175){\line(-1,0){90}}
\bezier45(-108,206.5)(-100,198)(-90,208)
\bezier20(-105,205)(-100,210)(-95,205)
\bezier45(-105,193)(-95,183)(-87,191.5)
\bezier20(-100,190)(-95,195)(-90,190)

\end{picture}
\vspace{1in}
\end{center}
    \caption{An asymptotically Delaunay surface $M \subset \R^3$.}%
    \label{fig:cmcsurfR3}
\end{figure}
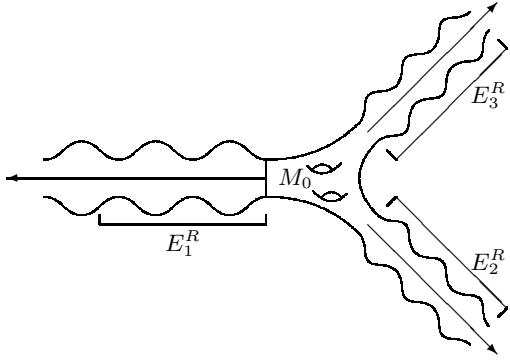

\paragraph{5.2 Asymptotically Delaunay hypersurfaces.}\label{SS:ads}
Let $M\subset \R^{\,n+1}$ be a hypersurface such that

\begin{enumerate}

\item $M$ can be decomposed as
$$M = M_0 \bigsqcup \, \bigsqcup_{j=1}^N E_j \; , $$
where $M_0$ is compact with boundary and where each $E_j$ is an
end of $M$ (Figure 3).

\item Each end $E_j$ is a cylindrical graph over half a Delaunay
  unduloid $\cD_j(\mu_j)$ with weight $\mu_j>0$ and semi-axis $a_j+
  \R_+\, d_j$ for some $a_j, d_j \in \R^{n+1}$. The boundary of $E_j$ lies
  in the hyperplane through $a_j$
  orthogonal to $d_j$, and $\partial M_0 = \bigsqcup_j \,\partial E_j$.

\item The graph $E_j$ above $\cD(\mu_j)$ is given by a
parametrization of the form
\begin{equation}\label{E:adh1a}
\R \times S^{n-1} \ni (x,\omega) \to F_j(x,\omega) = \big( x,
(f(x) + w_j(x,\omega))\,\omega \big)
\end{equation}
with some function $w_j(x,\omega)$, for $(x,\omega) \in \R_+\times
S^{n-1}$, where $f(x)$ satisfies equation (\ref{E:delau2a}).  We
assume that $w_j$ tends to zero in $C^2$-norm on $[r,\infty[\times
S^{n-1}$ as $r \to \infty$.
\end{enumerate}

\begin{definition}\label{D:ads}
We will say that a hypersurface which satisfies the preceding
three conditions is an {\em asymptotically Delaunay} hypersurface.
\end{definition}

\medskip

This definition extends mutatis mutandis to the case of
hypersurfaces in $\HH^{\,n+1}$ (in this case, the axis is a
geodesic ray parametrized by arc-length). Note that \cite{KKS} and
\cite{KKMS} give sufficient conditions to insure that a
cmc hypersurface is asymptotically Delaunay.

\medskip

With the above notations, we also introduce the following subsets
of $M$:
$$M^R = M_0 \bigsqcup \,
\bigsqcup_{j=1}^N E_j^R \; \;\;\;\;\;\;\;\;\; \mbox{for} \;\;\;
R>0 \; , $$ where $E_j^R$ is the part of $E_j$ which lies above
$a_j + [0,R]\,d_j$ (see Figure 3), and
$$M^{S,R} = M^R \setminus M^S = \bigsqcup_{j=1}^N E_j^{S,R} \;
\;\;\;\;\;\;\;\;\; \mbox{for} \;\;\; R>S>0 \; , $$
where $E_j^{S,R}$ is the part of $E_j$ which lies above $a_j +
[S,R]\,d_j$.
We can use similar notations for hypersurfaces $M$ in $\HH^{\,n+1}$.

\paragraph{5.3 Main results.}\label{SS:thms}
We have the following results:

\begin{theorem}\label{T:thm2}
Let $M \subset \R^{\, n+1}$ be a complete asymptotically
Delaunay cmc hypersurface. Let $E_j, j=1
\ldots N,$ be the ends of $M$ and let $\cD(\mu_j)$ be the Delaunay
unduloid to which $E_j$ is asymptotic (with weight parameter
$\mu_j> 0$). Denote by $T(\mu_j)$ the period of $\cD(\mu_j)$.  Then
\begin{equation}\label{E:ethm2}
\lim_{R\to \infty }\frac{\ind (M\cap B(R))}{R} = 2\,
\sum_{j=1}^{N}\frac{1}{T(\mu_j)} \; ,
\end{equation}
where $B(R)$ is the Euclidean ball of radius $R$ in $\R^{\,n+1}$.
\end{theorem}

\medskip

\textbf{Remark.~} It follows from \cite{KKS} that the preceding
theorem applies when $M$ is a properly embedded cmc $1$ surface with
finite topology in $\R^3$.

\medskip

Let $M \subset \HH^{\, n+1}$ be a complete properly embedded
hypersurface, with cmc $H > 1$ and finite topology. Such an $M$ is
asymptotically Delaunay if $n=3$ (\cite{KKMS}, Theorem 1.2) or if
$n\ge 4$ and each end of $M$ is within a bounded distance of a
geodesic ray (\cite{KKMS}, Theorem 1.3).

\begin{theorem}\label{T:thm3}
Let $M \subset \HH^{\, n+1}(-1)$ be a complete properly embedded
hypersurface, with cmc $H > 1$ and finite topology. Assume
furthermore that each end of $M$ is within a bounded distance of
some geodesic ray when $n\ge 4$. Let $E_j, j=1 \ldots N,$ be the
ends of $M$. Let $\cD_H(\mu_j)$ be the Delaunay unduloid (with
weight parameter $\mu_j > 0$) to which $E_j$ is asymptotic. Denote
by $\tau(\mu_j)$ the period of $\cD_H(\mu_j)$. Then
\begin{equation}\label{E:ethm3}
\lim_{R\to \infty }\frac{\ind (M\cap B(R))}{R} = 2 \,
\sum_{j=1}^{N}\frac{1}{\tau(\mu_j)} \; ,
\end{equation}
where $B(R)$ is the hyperbolic ball of hyperbolic radius $R$ in
$\HH^{\;n+1}$.
\end{theorem}

\medskip

\textbf{Remark.~} The general idea of the proofs of these theorems
is to apply \emph{Dirichlet--Neumann bracketing} (see \cite{RS} for example)
to $M^R$ decomposed as $M^R = M^S \sqcup M^{S,R}$,
and to use the fact that each component of $M^{S,R}$ (as $R \to \infty$)
is asymptotic to a Delaunay piece for which we can estimate the
index.

\medskip

\textbf{Proofs, main argument.~} Here we give the argument only for
the Euclidean
case.  (The argument in the hyperbolic case is identical, except for some
minor changes of notation.)  For $M$ as in Theorem \ref{T:thm2}, we
want to estimate the limits

$$ \liminf_{R \to \infty} \frac{\ind(M\cap B(R))}{R}=\liminf_{R \to \infty}
\frac{\ind(M^R)}{R}\;,
~~\limsup_{R \to \infty}
\frac{\ind(M\cap B(R))}{R}=\limsup_{R \to \infty} \frac{\ind(M^R)}{R} \; . $$

\medskip
\textbf{Step 1, Estimating the index from below.~}
Let $a_j$ be the scalar product of the normal to the hypersurface
with the Killing field corresponding to translation along the
axis of the Delaunay unduloid $\cD(\mu_j)$. Since the end $E_j$ is
asymptotic to $\cD(\mu_j)$, the nodal domains of the function $a_j$
look very much like the nodal domains of the corresponding
function for $\cD(\mu_j)$. Then, applying Dirichlet-Neumann bracketing to
the decomposition $M^R = M^S \sqcup M^{S,R}$ and letting $R \to \infty$,
Proposition \ref{P:ibp} implies $\liminf_{R \to \infty} \ind(M^R)/R$ is
greater than or equal to the value in the right hand side of (\ref{E:ethm2}).

\bigskip
\textbf{Step 2, Estimating the index from above.~}
For $R > S$, we decompose $M^R$ into pieces $M^R = M^S
\bigsqcup \, \bigsqcup_{j=1}^{N} E_{j}^{S,R}$, chosen in such a way
that the components of $\partial M^{S}$ lie above
boundaries of $\cC(\mu_j)$ basic Neumann blocks of the corresponding Delaunay
unduloids (this can be done with the correct choices of $M_0$ and $S$).

\smallskip

Fix some $\ell \in \N$.  Each piece $E_j^{S,R}$ can again be
decomposed into almost Delaunay pieces above $\cC_\ell(\mu_j)$
pieces of the Delaunay unduloids $\cD(\mu_j)$, plus a remainder
part. We write such a decomposition
$$E_j^{S,R} = \bigsqcup_{p=1}^{m_j}\widetilde{\cC}_{\ell,p}(\mu_j)
\, \bigsqcup \widetilde{\cR}_j \; , $$
where $\cR_j \subset \widetilde{\cC}_{\ell, m_j+1}(\mu_j)$.
Dirichlet--Neumann bracketing implies
$$ \ind_D(M^R) \le   \Nind_N(M^{S}) + \sum_{j=1}^{N}\left\aco
\sum_{p=1}^{m_j} \Nind_N(\widetilde{\cC}_{\ell,p}(\mu_j)) +
\NDind_{ND}(\widetilde{\cR}_j)\right\acf \; ,$$ where $\ind_D$
(resp. $\Nind_N$, $\NDind_{ND}$) stands for the index with
Dirichlet (resp. Neumann, mixed Neumann--Dirichlet) boundary
condition.

\medskip

The number $\ell$ being fixed, we can choose $S$ (and $R > S$) so
large that each piece $\widetilde{\cC}_{\ell,p}(\mu_j)$ is close
enough to a $\cC_{\ell}(\mu_j)$-piece so that
$\Nind_N(\widetilde{\cC}_{\ell,p}(\mu_j)) \le 2 \ell + c(n,H)$, by
Lemma \ref{L:l10} and Proposition \ref{P:inp}.

\medskip

We can now look at the extrinsic length $R$ and write, for each
end $E_j$,
$$ S + m_j \, \ell \, T(\mu_j) \le R \le  S + (m_j  +
1) \, \ell \, T(\mu_j) \; .$$
It follows that
$$
 \ind_D(M^R) \le  \Nind_N(M^{S}) +
\sum_{j=1}^{N} \left\{ (1 + \frac{c(n,H)}{2 \ell}) \,
\frac{2(R-S)}{T(\mu_j)} + \NDind_{ND}(\widetilde{\cC}_{\ell}(\mu_j))
\right\} $$
$$
\le  \Nind_N(M^{S}) +
\sum_{j=1}^{N} \left\{ (1 + \frac{c(n,H)}{2 \ell}) \,
\frac{2(R-S)}{T(\mu_j)} + \NDind_{ND}(\cC_{\ell}(\mu_j))+c(n,H)
\right\} \; . $$

Dividing the preceding inequality by $R$ and letting $R$ tend to
infinity, we find that

$$\limsup_{R\to \infty} \frac{\ind(M^{R})}{R} \le 2\, \sum_{j=1}^N (1 +
\frac{c(n,H)}{2 \ell}) \frac{1}{T(\mu_j)} \; . $$

\medskip

Since $\ell$ is an arbitrary positive integer, we have that
$\limsup_{R\to \infty} \ind(M^{R})/R$ is less than or equal to
the value in the right hand side of (\ref{E:ethm2}).


\paragraph{5.4 Other growth results.}\label{SS:ogr}
Kapouleas \cite{K} has constructed examples of complete constant mean
curvature surfaces in $\R^{3}$ which are periodic with respect to
some $2$ (resp. $3$) dimensional lattice.
It is not difficult to establish that,
for each of the doubly (resp. triply) periodic surfaces $M$ in
\cite{K}, there exist finite positive constants $c_1$ and $c_2$ such that
$c_1 R^2 \le \mbox{Ind}(M \cap B(R)) \le c_2 R^2$ (resp.
$c_1 R^3 \le \mbox{Ind}(M \cap B(R)) \le c_2 R^3$) for large $R$.

\vskip1cm

\parbox[t]{5.2in}{
Pierre Bérard\\ Institut Fourier, UMR 5582 UJF--CNRS, Université
Joseph Fourier\\ B.P. 74, 38402 St Martin d'Hères Cedex, France\\
Pierre.Berard@ujf-grenoble.fr\\[4pt]
\texttt{http://www-fourier.ujf-grenoble.fr/$\sim$pberard/publications.html}}

\bigskip

\parbox[t]{5.2in}{
Levi Lopes de Lima\\ Departamento de Matemática, Universidade
Federal do Ceará \\Campus do Pici, 60455--760 Fortaleza,
Brazil\\ levi@mat.ufc.br}

\bigskip

\parbox[t]{5.2in}{
Wayne Rossman\\ Department of Mathematics, Faculty of Science\\
Kobe University, Rokko, Kobe 657-8501, Japan\\
wayne@math.kobe-u.ac.jp\\[4pt]
\texttt{http://www.math.kobe-u.ac.jp/HOME/wayne/wayne.html}}

\end{document}